\tolerance=10000
\raggedbottom

\baselineskip=15pt
\parskip=1\jot

\def\sk{\vskip 3\jot}

\def\heading#1{\vskip3\jot{\noindent\bf #1}}
\def\label#1{{\noindent\it #1}}


\def\ref#1;#2;#3;#4;#5.{\item{[#1]} #2,#3,{\it #4},#5.}
\def\refinbook#1;#2;#3;#4;#5;#6.{\item{[#1]} #2, #3, #4, {\it #5},#6.} 
\def\refbook#1;#2;#3;#4.{\item{[#1]} #2,{\it #3},#4.}


\def\({\bigl(}
\def\){\bigr)}

\def\Ex{{\rm Ex}}
\def\Var{{\rm Var}}

\def\pwrt#1{{\partial\over\partial #1}}

\def\Li{{\rm Li}}

\def\ga{\gamma}
\def\ze{\zeta}
\def\la{\lambda}
\def\si{\sigma}
\def\ps{\psi}

\def\Ga{\Gamma}
\def\De{\Delta}

{
\pageno=0
\nopagenumbers
\rightline{\tt egp.ii.arxiv.tex}
\vskip1in

\centerline{\bf Asymptotic Behavior of the Moments of the Maximum Queue Length During a Busy Period}
\vskip0.5in

\centerline{Patrick Eschenfeldt}
\centerline{\tt peschenfeldt@hmc.edu}
\sk

\centerline{Ben Gross}
\centerline{\tt bgross@hmc.edu}
\sk

\centerline{Nicholas Pippenger}
\centerline{\tt njp@math.hmc.edu}
\sk

\centerline{Department of Mathematics}
\centerline{Harvey Mudd College}
\centerline{1250 Dartmouth Avenue}
\centerline{Claremont, CA 91711}
\vskip0.5in

\noindent{\bf Abstract:}
We give a simple derivation of the distribution of the maximum $L$ of the length of the queue during a busy period for the $M/M/1$ queue with $\la<1$ the ratio between arrival rate and service rate.
We observe that the asymptotic behavior of the moments of $L$ is related to that of Lambert series
for the generating functions for the sums of powers of divisors of positive integers.
We show that $\Ex[L]\sim \log\(1/(1-\la)\)$ and $\Ex[L^k]\sim k!\,\ze(k)/(1-\la)^{k-1}$ for $k\ge 2$, so that $\Var[L]\sim 3/\pi^2(1-\la)$.
More generally, we show how to obtain asymptotic expansions for these moments with error terms of the form $O\((1-\la)^N\)$ for any $N$.

\vskip0.5in
\leftline{{\bf Keywords:} Queueing theory, Lambert series, asymptotic expansions.}
\sk
\leftline{{\bf Subject Classification:} 60K26, 90B22}

\vfill\eject
}

\heading{1.  Introduction}

We study the $M/M/1$ queue (that is, the single server queue with independent exponentially distributed interarrival times and independent exponentially distributed service times).
Let $\la$ denote the ratio of the arrival rate (the reciprocal of the mean interarrival time) to the service rate (the reciprocal of the mean service time).
When $\la<1$, the the busy period (the interval during which the service is continuously busy)  is finite with probability one.
Let the random variable $L$ denote the maximum length of the queue during a busy period.
(We count the customer being served in the queue length, so that $L\ge 1$ and $L=1$ when and only when the busy period comprises a single service interval.)
It is known that $L$ has the distribution
$$\Pr[L>l] = {(1-\la)\,\la^l \over 1-\la^{l+1}}; \eqno(1.1)$$
see for example Cohen [C1, pp.~191--193].
(Neuts [N] has also treated this question, but does not give a simple formula for the distribution.)
In Section 2, we shall give a simple derivation of (1.1) based on the solution to the ``Gambler's Ruin''  problem.

Our main goal in this paper is to give asymptotic expansions for the moments 
$$\Ex[L^k] = \sum_{l\ge 0} l^k \, \Pr[L=l]$$
of $L$ in the ``heavy traffic'' limit as $\la\to 1$.
Writing 
$$\eqalign{
\De_k(m) 
&= m^k - (m-1)^k \cr
&= \sum_{0\le j\le k-1} {k\choose j}(-1)^{k-1-j} \, m^j \cr
}$$
for the backward differences of the $k$-th powers, and setting
$$S_k(\la) = \sum_{m\ge 1} {m^k \la^m \over 1 - \la^m}, \eqno(1.2)$$
summation by parts yields
$$\eqalignno{
\Ex[L^k]
&= \sum_{l\ge 0} l^k \, \Pr[L=l] \cr
&= \sum_{l\ge 0} \De_k(l+1) \,\Pr[L>l] \cr
&= (1-\la) \sum_{l\ge 0} {\De_k(l+1) \,\la^l \over 1-\la^{l+1}} \cr
&= {1-\la\over \la} \sum_{m\ge 1} {\De_k(m) \,\la^m \over 1-\la^m} \cr
&= {1-\la\over \la} \sum_{m\ge 1} \sum_{0\le j\le k-1} {k\choose j} (-1)^{k-1-j} {m^j \,\la^m \over 1-\la^m} \cr
&= {1-\la\over \la} \sum_{0\le j\le k-1} {k\choose j} (-1)^{k-1-j} S_j(\la). &(1.3)\cr
}$$
Since $\Ex[L^k]$ is a linear combination of the $S_k(\la)$, it will suffice to determine the asymptotic behavior of the sums $S_k(\la)$.

The sums in (1.2), which are called Lambert series, arise in a natural way in number theory
(see for example Hardy and Wright [H, p.~257]).
We have
$$\eqalignno{
S_k(\la)
&=\sum_{m\ge 1} {m^k \, \la^m \over 1-\la^m} \cr
&= \sum_{m\ge 1} m^k \, \sum_{l\ge 1} \la^{lm} \cr
&= \sum_{l\ge 1} \sum_{m\ge 1} m^k \, \la^{lm} \cr
&= \sum_{n\ge 1} \,\la^n \, \sum_{d\mid n} d^k &(1.4)\cr
&= \sum_{n\ge 1} \si_k(n)\,\la^n, \cr
}$$
where the inner sum in (1.4) is over integers $d$ dividing $n$, and
$\si_k(n)$ denotes the sum of the $k$-th powers of the divisors of $n$
(see Hardy and Wright [H, p.~239]).
Thus $S_k(\la)$ is the generating function for $\si_k(n)$.

We note that the Lambert series $S_k(\la)$ can be expressed in terms of known (albeit exotic) functions of analysis.
We define the $q$-gamma function by
$$\Ga_q(x) = (1-q)^{1-x} \prod_{n\ge 0} {1-q^{n+1}\over 1-q^{n+x}}$$
(see for example Gasper and Rahman[G, p.~16]). 
(This function gets its name from the fact that
$\lim_{q\to 1} \Ga_q(x) = \Ga(x)$, where $\Ga(x)$ is the Euler gamma function; see for example
Whittaker and Watson [W, pp.~235--264].)
If we define the $q$-digamma function $\ps_q(x)$ as the logarithmic derivative
$$\eqalign{
\ps_q(x)
&= \pwrt{x} \log \Ga_q(x) \cr
&= -\log(1-q) + \log q \, \sum_{n\ge 0} {q^{n+x} \over 1-q^{n+x}} \cr
}$$
of the $q$-gamma function, then we have
$$S_0(\la) = {\ps_\la(1) + \log(1-\la) \over \log\la}.$$
To go further, we define the $k$-th $q$-polygamma function $\ps_q^{(k)}$ as the $k$-th derivative
$$\ps_q^{(k)}(x) = \left(\pwrt{x}\right)^k \ps_q(x)$$
of the $q$-digamma function.
If we set $z=q^{n+x}$, then
$$\left(z\pwrt{z}\right) = \left({1\over \log q}\,\pwrt{x}\right).$$
Since
$$\sum_{m\ge 1} m^k z^m = \left(z\pwrt{z}\right)^k {z\over 1-z},$$
we have
$$\sum_{m\ge 1} m^k \, q^{(n+x)m} 
= {1\over \log^k q} \,  \left(\pwrt{x}\right)^k \, {q^{n+x} \over 1 - q^{n+x}}.$$
Summing over $n\ge 0$ yields
$$\eqalign{
\sum_{m\ge 1} {m^k \, q^{xm} \over 1 - q^{xm}}
&= \sum_{m\ge 1} m^k \,  \sum_{n\ge 0}  q^{m(n+x)} \cr
&= {1\over \log^k q} \, \left(\pwrt{x}\right)^k \, 
\sum_{n\ge 0} {q^{n+x} \over 1 - q^{n+x}} \cr
&= {1\over \log^k q} \, \left(\pwrt{x}\right)^k \, 
 {\ps_q(x) + \log(1-q) \over \log q}. \cr
}$$
Thus for $k\ge 1$ we have
$$S_k(\la) = {\ps_\la^{(k)}(1) \over \log^{k+1} \la}.$$

In Section 3, we shall begin our study of the asymptotic behavior of the moments of $L$,
deriving the leading terms
$$\Ex[L] \sim \log {1\over 1-\la}, \eqno(1.5)$$
and, for $k\ge 2$,
$$\Ex[L^k] \sim {k!\,\ze(k) \over (1-\la)^{-1}}, \eqno(1.6)$$
where $\ze(k) = \sum_{n\ge 1} 1/n^k$ is the Riemann zeta function.
It will be noted that $\Ex[L]$ grows quite slowly as $\la\to 1$.
If the random variable $K$ denotes the length of the queue in equilibrium, then 
$\Ex[K] = \la/(1-\la)$, which grows much more rapidly
(see Cohen [C], p.~181).
It may appear paradoxical that the maximum queue length grows more slowly than
the mean queue length, but it must be borne in mind that $\Ex[L]$ is an average over busy 
periods, whereas $\Ex[K]$ is an average over time.
Indeed, the majority of busy periods have $L=1$: after the arrival initiating the busy period, 
the next event determines whether $L=1$ (if that event is a service termination) or $L>1$
(if that event is another arrival).
Because  $\la<1$,  the former (with probability $1/(1-\la)$) is more likely than the latter (with probability $\la/(1-\la)$).
We also note that, since $\ze(2) = \pi^2/6$, we have $\Var[L] = \Ex[L^2] - \Ex[L]^2 \sim \pi^2/3(1-\la)$, which grows much more rapidly 
than $\Ex[L]$ (or even $\Ex[L]^2$).

In Section 4, we shall refine these results by showing how to obtain complete asymptotic expansions for the moments $\Ex[L^k]$ as $\la\to 1$.
We obtain error terms of the form $O\((1-\la)^c\)$ for any $c$.
These refinements are obtained by methods recently introduced in rigorous quantum field theory.
All our results can be extended to the $M/M/s$ queue (with $s$ independent and identical servers), if the busy period is defined as a contiguous interval during which all $s$ servers are busy, and 
$\la$ is the ratio of the arrival rate to $s$ times the service rate for each server; for simplicity we confine our attention to the case $s=1$.
\vfill\eject

\heading{2. The Distribution}

In this section, we shall give a simple derivation of the formula (1.1).
Consider a game played between two players: $P$, who begins with $v$ dollars, and $Q$ who begins with $w$ dollars.
At each step of the game, a biassed coin is tossed; $P$ wins with probability $p$, in which case 
$Q$ pays $P$ one dollar, and $Q$ wins with the complementary probability $q=1-p$, in which case 
$P$ pays $Q$ one dollar.
The game continues until one of the players is ruined (that is, has no money left).
It is known that (1) with probability one, either $P$ or $Q$ is eventually ruined, and (2),
if $p\not=q$, then the probability that $Q$ is ruined is
$$\Pr[Q\hbox{\ ruined}] = {(q/p)^v - 1 \over (q/p)^{v+w} - 1} \eqno(2.1)$$
(see for example Feller [F, p.~345]).

Now consider a busy period of the $M/M/1$ queue.
The successive events of arrivals and terminations of service intervals during the busy period correspond to steps in the game described above. 
The wealth of player $P$ will correspond to the length of the queue at each step, so $v=1$.
An arrival will correspond to a win by player $P$, so $p=\la/(1+\la)$, and
the termination of a service interval will correspond to a win by player $Q$, so $q=1/(1+\la)$.
Suppose that player $Q$ begins with $w=l$ dollars.
Then the event $L>l$ will correspond to $Q$ being ruined.
Substituting these values in (2.1) yields (1.1).

This correspondence also shows what happens for $\la\ge 1$.
For $\la=1$ (in which case the busy period is finite with probability one, but its expected length is infinite), we have take $p=q=1/2$, and have
$$\Pr[Q\hbox{\ ruined}] = {v\over v+w}.$$
This result yields
$$\Pr[L>l] = {1\over l+1},$$
so that 
$$\Ex[L] = \sum_{l\ge 0} \Pr[L>l] \eqno(2.2)$$ 
diverges logarithmically.
Of course, for $\la>1$ (in which case the busy period is infinite with positive probability),
(2.1) shows that (2.2) diverges linearly.
\sk

\heading{3. The Leading Terms of the Moments}

In this section, we shall derive the leading terms (1.5) and (1.6) of the asymptotic expansions for the moments of $L$.
We begin by deriving the leading terms of the asymptotic expansions for the sums
$$S_0(\la) \sim {1\over 1-\la}\log {1\over 1-\la} \eqno(3.1)$$
and, for $j\ge 1$,
$$S_j(\la) \sim {j! \, \ze(j+1)\over (1-\la)^{j+1}}. \eqno(3.2)$$
Once these formulas are established, it will be clear that the sum in (1.3) is dominated by the term for which $j=k-1$, so that $\Ex[L^k] \sim k\,S_{k-1}(\la)$, and (1.5) and (1.6) follow from (3.1) and (3.2), respectively.

Our strategy for proving (3.1) and (3.2) will be to approximate the sums
$S_j(\la)$
by integrals
$$I_j(\la) = \int_1^\infty {x^j \, \la^x \, dx \over 1-\la^x },$$
then then to show that the difference $S_j(\la)-I_j(\la)$ is negligible in comparison with $I_j(\la)$.
It will be convenient to write $\la = e^{-h}$.
The limit $\la\to 1$ then corresponds to $h\to 0$.
We have
$$\eqalignno{
h
&= \log{1\over \la} \cr 
&= \log{1\over 1-(1-\la)} \cr
&\sim 1-\la. &(3.3) \cr
}$$

For $j=0$, we have
$$\eqalign{
I_0(\la) 
&= \int_1^\infty { \la^x \, dx \over 1-\la^x } \cr
&= \int_1^\infty  \sum_{l\ge 1}  \, \la^{lx} \,dx  \cr
&=  \sum_{l\ge 1}  \, \int_1^\infty  e^{-hlx} \,dx  \cr
&=  \sum_{l\ge 1}  \, {e^{-hl} \over hl }  \cr
&=  {1\over h} \sum_{l\ge 1}  \, {\la^{l} \over l }  \cr
&= {1 \over h }\;\log {1\over 1-\la}. \cr
}$$
Substituting (3.3) in this result yields 
$$I_0(\la) \sim  {1\over 1-\la}\log {1\over 1-\la}. \eqno(3.4)$$
We bound $\vert S_0(\la) - I_0(\la)\vert$ by the total variation of $f(x) = \la^x / (1-\la^x)$.
Since $f(x)$ decreases monotonically from $\la/(1-\la)$ to $0$ as $x$ increases from 
$1$ to $\infty$, we have $\vert S_0(\la) - I_0(\la)\vert \le \la/(1-\la) \sim 1/(1-\la)$.
Since this difference is negligible in comparison with (3.4), we obtain (3.1).

For $j\ge 1$, we have
$$\eqalignno{
I_j(\la)
&= \int_1^\infty {x^j \, \la^x \, dx \over 1-\la^x} \cr
&= \int_0^\infty {(y+1)^j \, \la^{y+1} \, dy \over 1-\la^{y+1}} \cr
&= \int_0^\infty \sum_{0\le i\le j} {j\choose i}{y^i \, \la^{y+1} \, dy \over 1-\la^{y+1}} \cr
&= \int_0^\infty \sum_{0\le i\le j} {j\choose i}y^i \sum_{l\ge 1} \la^{l(y+1)} \, dy  \cr
&= \int_0^\infty \sum_{0\le i\le j} {j\choose i}y^i \sum_{l\ge 1} e^{-hl(y+1)} \, dy  \cr
&=  \sum_{0\le i\le j} {j\choose i}\sum_{l\ge 1}\int_0^\infty y^i \, e^{-hl(y+1)} \, dy  \cr
&=  \sum_{0\le i\le j} {j\choose i}\sum_{l\ge 1} {i! \, e^{-hl} \over (hl)^{i+1}} \cr
&= \sum_{0\le i\le j} {j\choose i} {i! \over h^{i+1}} \sum_{l\ge 1} {\la^{l} \over l^{i+1}} \cr
&= \sum_{0\le i\le j} {j\choose i} {i! \over h^{i+1}} \Li_{i+1}(\la), &(3.5)\cr
}$$
where $\Li_k(\la) = \sum_{n\ge 1} \la^n/n^k$ is the $k$-th polylogarithm.
Since $\Li_1(\la) = \log \(1/(1-\la)\)$ and $\Li_k(\la)\to\ze(k)$ as $\la\to 1$ for $k\ge 2$,
the sum in (3.5) is dominated by the term for which $i=j$, and we have
$$\eqalignno{
I_j(\la) 
&\sim {j! \, \ze(j+1) \over h^{j+1}} \cr
&\sim {j! \, \ze(j+1) \over (1-\la)^{j+1}}&(3.6) \cr
}$$
We bound $\vert S_j(\la)-I_j(\la)\vert$ by the total variation of $f(x) = x^j \, \la^x / (1-\la^x)$ for
$0\le x < \infty$.
As $x$ increases, $f(x)$ increases monotonically from $0$ to a maximum, then decreases
monotonically to $0$.
Thus the total variation of $f(x)$ is twice the maximum.
This maximum is
$$\eqalign{
\max_{0\le x < \infty} f(x)
&= \max_{0\le x < \infty} {x^j \, e^{-hx} \over 1 - e^{-hx}} \cr
&= \max_{0\le x < \infty} {x^j  \over e^{hx} - 1} \cr
&= {1\over h^j} \, \max_{0\le y < \infty}  {y^j \over e^y - 1}. \cr
}$$
Furthermore, $y^j / (e^y - 1) \le j!$, because $e^y - 1 = \sum_{n\ge 1} y^n/n! \ge y^j/j!$.
Thus $\vert S_j(\la)-I_j(\la)\vert \le 2\max_{0\le x<\infty} f(x) \le 2j!/h^j \sim 2j!/(1-\la)^j$.
Since this difference is negligible in comparison with (3.6), we obtain (3.2).
\sk

\heading{4. The Complete Asymptotic Expansions}

In this section we shall show how asymptotic expansions, with error terms of the form
$O\((1-\la)^N\)$ for any $N$, can be derived for all of the moments $\Ex[L^k]$.
The essence of the argument is to use the Euler-Maclaurin formula to estimate the difference
between $S_j(\la)$ and $I_j(\la)$.
This is most conveniently done using a result of Zagier [Z].
Indeed, for $j\ge 1$, Zagier gives the expansion for $S_j(\la)$, in terms of the parameter
$h = -\log\la$ rather than $1-\la$.
All that remains for us to do is substitute an expansion for $h$ in terms of $1-\la$.
For $j=0$, the expansion for $S_0(\la)$ in terms of $h$ has been given by Egger (n\'{e} Endres) and Steiner [E1, E2], again using the result of Zagier.
We shall proceed differently, to obtain an expansion involving $-\log(1-\la)$ rather than $-\log h$.

\label{Proposition:}
(Zagier [Z, p.~318])
Let $f(x)$ be analytic at $x=0$,
with power series $f(x) = \sum_{n\ge 0} b_n \, x^n$ about $x=0$.
Suppose that $\int_0^\infty \vert f^{(N)}(x)\vert \, dx < \infty$ for all $N\ge 0$,
where $f^{(N)}(x)$ denotes the $N$-th derivative of $f(x)$.
Define $I_f = \int_0^\infty f(x) \, dx$.
Let $g(x) = \sum_{m\ge 1} f(mx)$.
Then $g(x)$ has the asymptotic expansion
$$g(x) \sim {I_f \over x} + \sum_{n\ge 0} {b_n \, B_{n+1} \, (-1)^n x^n \over (n+1) }, \eqno(4.1)$$
where $B_k$ is the $k$-th Bernoulli number, defined by $t/(e^t - 1) = \sum_{k\ge 0} B_k \, t^k/k!$.

This result is proved by using the Euler-Maclauren formula,
$$\eqalign{
\int_0^M f(y)\,dy
&= {f(0)\over 2} + \sum_{1\le m\le M-1} f(m) + {f(M)\over 2}
+ \sum_{1\le n\le N-1} {(-1)^n B_{n+1} \over (n+1)!}\(f^{(n)}(M) - f^{(n)}(0)\) \cr
&\qquad +(-1)^N \int_0^M f^{(N)}y {{B}_N(\{y\})\over N!}\,dy, \cr
}$$
where ${B}_k(y)$ is the  $k$-th Bernoulli polynomial, defined by 
$t e^{yt}/(e^t - 1) = \sum_{k\ge 0} B_k(y) \, t^k/k!$, and $\{y\} = y-\lfloor y\rfloor$ denotes the fractional part of $y$ .
(For the Euler-Maclauren formula, the Bernoulli numbers and the Bernoulli polynomials,
see for example Whittaker and Watson [W, pp.~125--128], where, however, the indexing of the numbers and polynomials is different.)
The condition $\int_0^\infty \vert f^{(N)}(y)\vert \, dy < \infty$ allows us to let $M\to\infty$, obtaining
$$\eqalign{
\int_0^\infty f(y)\,dy 
&=  \sum_{m\ge 1} f(m) 
+ \sum_{0\le n\le N-1} {(-1)^n B_{n+1} \over (n+1)!} \, f^{(n)}(0) 
 +(-1)^N \int_0^\infty f^{(N)}(y) {{B}_N(\{y\})\over N!}\,dy. \cr
}$$
If we now write $f(xy)$ instead of $f(y)$, we obtain
$$\eqalign{
\int_0^\infty f(xy)\,dy 
&=  \sum_{m\ge 1} f(mx) 
+ \sum_{0\le n\le N-1} {(-1)^n B_{n+1} \over (n+1)!} \, f^{(n)}(0) \, x^n
 +(-1)^N \, x^N \int_0^\infty f^{(N)}(xy) {{B}_N(\{y\})\over N!}\,dy. \cr
}$$
Changing the variable of integration from $y$ to $y/x$ then yields
$$\eqalign{
{1\over x} \int_0^\infty f(y)\,dy 
&=  \sum_{m\ge 1} f(mx) 
+ \sum_{0\le n\le N-1} {(-1)^n B_{n+1} \over (n+1)!} \, f^{(n)}(0) \, x^n
 +(-1)^N \, x^{N-1} \int_0^\infty f^{(N)}(y) {{B}_N(\{y/x\})\over N!}\,dy. \cr
}$$
The integral on the left-hand side id $I_f$, the first sum on the right-hand side is $g(x)$,
$f^{(n)}(0) = n!\,b_n$, and the last term on the right-hand side is $O(x^{N-1})$.
Thus
$${I_f\over x} = g(x) + \sum_{0\le n\le N-1} {b_n \, B_{n+1} \, (-1)^n x^n \over (n+1) }
+ O(x^{N-1}),$$
which yields the expansion (4.1).

For $j\ge 1$, we define
$$f(x) = {x^j \over e^x -1}.$$
Then $f(x)$ is analytic at $x=0$ with the Taylor series
$$f(x) = \sum_{n\ge 0} {B_n \, x^{n+j-1} \over n!}$$
and the integral
$$\eqalign{
I_f
&= \int_0^\infty {x^j \, e^{-x} \, dx \over 1 - e^{-x}} \cr
&= j! \, \ze(j+1) \cr
}$$
(see for example Whittaker and Watson [W, p.~266]).
Furthermore, $f^{(N)}(x)$ is a rational function of $x$ and $e^x$, in which the degree of the numerator in $e^x$ is $N$, while the denominator is $(e^x-1)^{N+1}$.
Thus $f(x)$ satisfies the conditions of the proposition, and
we have the asymptotic expansion
$$g(x) \sim {j! \, \ze(j+1) \over x} + 
\sum_{n\ge 0} {(-1)^{n+j-1} \, B_n \, B_{n+j} \, x^{n+j-1} \over n! \, (n+j)}.$$
Recalling that $\la = e^{-h}$, so that $h=-\log\la$, we therefore have
$$\eqalignno{
S_j(\la)
&= \sum_{m\ge 1} {m^j \, e^{-hm} \over 1 - e^{-hm}} \cr
&= {1\over h^j} \sum_{m\ge 1} {(mh)^j  \over e^{hm} - 1} \cr
&= {1\over h^j} \sum_{m\ge 1} f(mh) \cr
&= {g(h) \over h^j}  \cr
&\sim {j! \, \ze(j+1) \over h^{j+1}} + 
\sum_{n\ge 0} {(-1)^{n+j-1} \, B_n \, B_{n+j} \, h^{n-1} \over n! \, (n+j)}. &(4.2)\cr
}$$
We note that, if $j$ is odd, then this expansion has only finitely many terms (because $B_n=0$ for odd $n\ge 3$).
To obtain an asymptotic expansion in terms of $1-\la$, we must substitute the expansion for $1/h$:
$$\eqalignno{
{1\over h}
&= {1\over -\log \la} \cr
&= {-1\over \log \(1-(1-\la)\)} \cr
&= {1\over 1-\la}\, {-(1-\la)\over \log \(1-(1-\la)\)} \cr
&= {1\over 1-\la}\, \sum_{n\ge 0} {(-1)^n \, C_n \, (1-\la)^n \over n!}, &(4.3)\cr
}$$
where $C_k$ is the $k$-th Bernoulli number of the second kind, defined by
$t/\log(1+t) = \sum_{k\ge 0} C_k \, t^k/k!$ (see for example Roman [R, p.~116]).
(These numbers are also called the Cauchy numbers of the first kind, and are given by
$C_k = \int_0^1 x(x-1)\cdots(x-k+1)\,dx$; see for example Comtet [C2, pp.~293--294].)

For $j=0$, we must proceed differently, because 
$$f(x) = {1\over e^x - 1}$$
has a pole at $x=0$.
We define
$$\eqalign{
f^*(x) 
&= f(x)  - {e^{-x} \over x} \cr
&= {1\over e^x - 1} - {e^{-x} \over x}. \cr
}$$
Then $f^*(x)$ is analytic at $x=0$ with the Taylor series
$$f^*(x) = \sum_{n\ge 0} {\(B_{n+1} - (-1)^{n+1}\) \, x^n \over (n+1)!}$$
and the integral
$$\eqalign{
I_{f^*}
&= \int_0^\infty  \left({1\over e^x - 1} - {e^{-x} \over x}\right) \, dx\cr
& = \ga \cr
}$$
(see for example Whittaker and Watson [W, p.~246]).
Furthermore, $f^{*(N)}(x)$ is a rational function of $x$ and $e^x$, in which the degree of the numerator in $e^x$ is $N$, while the denominator is $\((e^x-1)\,x\)^{N+1}$.
Thus $f^*(x)$ satisfies the conditions of the proposition, and
we have the asymptotic expansion
$$g^*(x) \sim {\ga\over x} 
+ \sum_{n\ge 0} {(-1)^n \, B_{n+1} \(B_{n+1} - (-1)^{n+1}\) \, x^n \over (n+1)\,(n+1)!}.$$
We therefore have
$$\eqalignno{
S_0(\la)
&= \sum_{m\ge 1} {e^{-mh} \over 1 - e^{-mh}} \cr
&= \sum_{m\ge 1} {1 \over e^{mh} - 1} \cr
&= \sum_{m\ge 1} {e^{-mh} \over mh} + \sum_{m\ge 1} {1 \over e^{mh} - 1} - {e^{-mh} \over mh} \cr
&= {1\over h}\log {1\over 1-\la} + \sum_{m\ge 1} f^*(mh) \cr
&= {1\over h}\log {1\over 1-\la} +  g^*(h) \cr
&\sim {1\over h}\log {1\over 1-\la} + {\ga\over h} 
+ \sum_{n\ge 0} {(-1)^n \, B_{n+1} \(B_{n+1} - (-1)^{n+1}\) \, h^n \over (n+1)\,(n+1)!}. &(4.4)\cr
}$$

To obtain asymptotic expansions for the moments of $L$, we substitute (4.3) into (4.2) and (4.4),
then substitute the results into (1.3),
using the expansion
$$\eqalign{
{1-\la\over\la}
&= {1-\la\over 1-(1-\la)} \cr
&= \sum_{n\ge 1} (1-\la)^n. \cr
}$$
Retaining only terms that do not vanish as $\la\to 1$, we obtain
$$\Ex[L] = \log{1\over 1-\la} + \ga + O\left((1-\la)\log{1\over 1-\la}\right)$$
and
$$\Ex[L^2] = {3\pi^2\over 6(1-\la)} + \log{1\over 1-\la} 
+ (\ga-1) + O\left((1-\la)\log{1\over 1-\la}\right)$$
for the first two moments.
Thus we have
$$\eqalign{
\Var[L]
&= \Ex[L^2] - \Ex[L]^2 \cr
&= {3\pi^2\over 6(1-\la)} - \log^2 {1\over 1-\la} + (1-2\ga)\log{1\over 1-\la} - \ga^2
+O\left((1-\la)\log^2 {1\over 1-\la}\right). \cr
}$$
\sk

\heading{5. Acknowledgment}

The research reported here was supported
by Grant CCF  0917026 from the National Science Foundation.
\sk

\heading{6. References}

\refbook C1;  J. W. Cohen;
The Single Server Queue (revised edition);
North-Holland, Amsterdam, 1982.

\refbook C2; L. Comtet;
Advanced Combinatorics:
The Art of Finite and Infinite Expansion;
D.~Reidel Publishing Co., Dortrecht, 1974. 

\ref E1; S. Egger (n\'{e} Endres) and F. Steiner;
``A New Proof of the Vorono\"{\i} Summation Formula'';
J. Phys.\ A: Math.\ Theor.; 44 (2011) 225302 (11 pp.).

\ref E2; S. Endres and F. Steiner;
``A Simple Infinite Quantum Graph'';
Ulmer Seminare Funktionalanalysis und Differentialgleichungen'';
14 (2009) 187--200.

\refbook F; W. Feller;
An Introduction to Probability Theory and Its Applications (3rd edition);
John Wiley \& Sons, New York, 1968.

\refbook G; G. Gasper and M. Rahman;
Basic Hypergeometric Series;
Cambridge University Press, Cambridge, 1990.

\refbook H; G. H. Hardy and E. M. Wright;
Introduction to the Theory of Numbers (5th edition);
Clarendon Press, Oxford, 1979.

\ref N; M. F. Neuts;
``The Distribution of the Maximum Length of  Poisson Queue During a Busy Period'';
Oper.\ Res.; 12:2 (1964) 281--285.

\refbook R; S. Roman;
The Umbral Calculus;
Academic Press, New York, 1984.

\refbook W; E. T. Whittaker and G. N. Watson;
A Course of Modern Analysis (4th edition);
Cambridge University Press, London, 1927.

\refinbook Z; D. Zagier;
``The Mellin Transformation and Other Useful Analytic Techniques'';
in: E.~Zeidler; Quantum Field Theory: I, Basics in Mathematics and Physics;
Springer-Verlag, Berlin, 2006, pp.~305--323.

\bye